\documentstyle[12pt]{article}

\def\o{\overline}

\makeatletter
\long\def\UR#1{\leavevmode\setbox\@tempboxa\hbox{#1}\@tempdima\fboxrule
    \advance\@tempdima \fboxsep \advance\@tempdima \dp\@tempboxa
   \hbox{\lower \@tempdima\hbox
  {\vbox{\hrule \@height \fboxrule
          \hbox{  \hskip\fboxsep
          \vbox{\vskip\fboxsep \box\@tempboxa\vskip\fboxsep}\hskip
                 \fboxsep\vrule \@width \fboxrule}%
                  }}}}
\makeatother

\makeatletter
\long\def\LR#1{\leavevmode\setbox\@tempboxa\hbox{#1}\@tempdima\fboxrule
    \advance\@tempdima \fboxsep \advance\@tempdima \dp\@tempboxa
   \hbox{\lower \@tempdima\hbox
  {\vbox{ 
          \hbox{  \hskip\fboxsep
          \vbox{\vskip\fboxsep \box\@tempboxa\vskip\fboxsep}\hskip
                 \fboxsep\vrule \@width \fboxrule}%
                 \hrule \@height \fboxrule}}}}
\makeatother

\makeatletter
\long\def\UL#1{\leavevmode\setbox\@tempboxa\hbox{#1}\@tempdima\fboxrule
    \advance\@tempdima \fboxsep \advance\@tempdima \dp\@tempboxa
   \hbox{\lower \@tempdima\hbox
  {\vbox{\hrule \@height \fboxrule
          \hbox{\vrule \@width \fboxrule \hskip\fboxsep
          \vbox{\vskip\fboxsep \box\@tempboxa\vskip\fboxsep}\hskip
                 \fboxsep }%
                  }}}}
\makeatother

\makeatletter
\long\def\LL#1{\leavevmode\setbox\@tempboxa\hbox{#1}\@tempdima\fboxrule
    \advance\@tempdima \fboxsep \advance\@tempdima \dp\@tempboxa
   \hbox{\lower \@tempdima\hbox
  {\vbox{ 
          \hbox{\vrule \@width \fboxrule \hskip\fboxsep
          \vbox{\vskip\fboxsep \box\@tempboxa\vskip\fboxsep}\hskip
                 \fboxsep }%
                 \hrule \@height \fboxrule}}}}
\makeatother

\let \ttorg \tt \def \tt{\ttorg \obeyspaces}

\begin{document}

\date{}

\title{\Large\bf Bi-oriented Quantum Algebras, and a Generalized Alexander Polynomial for Virtual Links}
\author{Louis Kauffman\thanks{Research supported in part by NSF Grant
DMS 980 2859} $\;\;$
and $\;\;$
David E. Radford\thanks{Research supported in part by NSF Grant
DMS 980 2178}\\ 
Department of Mathematics, Statistics \\
and Computer Science (m/c 249)    \\
851 South Morgan Street   \\
University of Illinois at Chicago\\
Chicago, Illinois 60607-7045}


\maketitle

\thispagestyle{empty}

\subsection*{\centering Abstract}

{\em This paper discusses the construction of a generalized Alexander polynomial for virtual knots and links, and 
the reformulation of this invariant as a quantum link invariant. We then introduce the concept of a bi-oriented 
quantum algebra which provides an algebraic context for this structure.}

\section{Introduction}
In this paper we discuss the construction of a generalized Alexander polynomial $G_{K}(s,t)$ via the concept of a biquandle.
Our approach leads directly to a generalization of the Burau representation upon which this invariant is based.
We then reformulate the invariant as a quantum link invariant and as a state summation.  In this context we show that 
the normalized quantum invariant $Z_{K}(\sigma, \tau)$ satisfies a Conway skein identity and reproduces $G_{K}(s,t)$.  
These invariants are useful for the theory of 
virtual knots and links as they vanish on classical knots and links.  Hence the invariants studied here can be used to show 
that many virtual knots and links are not classical. We give such examples in Section $2$. We also give an example of a virtual knot
that cannot be detected by
the generalized Alexander polynomial. This knot is detected by the structure of the corresponding generalized Alexander module.
We conclude with a diagram due to Kishino that is conjectured to be knotted, but so far has not been detected by any known invariants of
virtual knots. 

In the final section of the paper, we formulate the concept of a bi-oriented quantum
algebra.  This generalizes our previous concept of oriented quantum algebra \cite{KRO, KRCAT}  to include the necessary
structures to create invariants of virtual links. The invariant $Z_{K}(\sigma, \tau)$ studied in this paper fits non-trivially
into this framework.  
Subsequent papers will study the structure and applications of bi-oriented quantum algebras. 
\vspace{3mm}

\section{Virtual Links and a Generalized Alexander Polynomial}
In this section we will construct a generalization of the Alexander module and Alexander polynomial 
that is defined for virtual knots and links \cite{VKT, SVKT, DVK}. We then show how this generalized Alexander
polynomial can be seen as state summation model using a solution to the Quantum Yang-Baxter equation.
This state model will form the basis for the rest of the paper.
\bigbreak

Recall that classical knot theory can be described in terms of knot and link diagrams.
A {\em diagram} is a 4-regular plane graph (with extra structure at its nodes representing
the crossings in the link) represented on the surface of a plane and implicitly on the surface of
a two-dimensional sphere $S^{2}$.
One says that two such diagrams are {\em equivalent} if there is a sequence of moves of the
types indicated in part (A) of Figure 1 (The Reidemeister Moves) taking one diagram to the other.
These moves are performed locally on the 4-regular plane graph (with extra structure) that 
constitutes the link diagram.
\bigbreak

{\tt    \setlength{\unitlength}{0.92pt}
\begin{picture}(465,459)
\thicklines   \put(247,84){\vector(1,0){19}}
              \put(261,84){\vector(-1,0){17}}
              \put(348,199){\vector(1,0){19}}
              \put(362,199){\vector(-1,0){17}}
              \put(334,321){\vector(1,0){19}}
              \put(348,321){\vector(-1,0){17}}
              \put(334,419){\vector(1,0){19}}
              \put(348,419){\vector(-1,0){17}}
              \put(105,208){\vector(1,0){19}}
              \put(119,208){\vector(-1,0){17}}
              \put(91,322){\vector(1,0){19}}
              \put(105,322){\vector(-1,0){17}}
              \put(91,420){\vector(1,0){19}}
              \put(105,420){\vector(-1,0){17}}
              \put(239,1){\makebox(24,25){$C$}}
              \put(413,118){\makebox(23,26){$B$}}
              \put(42,127){\makebox(20,21){$A$}}
              \put(123,29){\framebox(262,99){}}
              \put(240,148){\framebox(224,310){}}
              \put(1,151){\framebox(226,307){}}
              \put(353,59){\circle{18}}
              \put(312,60){\circle{18}}
              \put(192,98){\circle{18}}
              \put(153,98){\circle{18}}
              \put(332,40){\line(1,1){39}}
              \put(291,81){\line(1,-1){41}}
              \put(171,120){\line(1,-1){40}}
              \put(133,80){\line(1,1){38}}
              \put(336,89){\line(1,1){34}}
              \put(292,41){\line(1,1){35}}
              \put(292,120){\line(1,-1){79}}
              \put(132,120){\line(1,-1){79}}
              \put(132,41){\line(1,1){35}}
              \put(177,86){\line(1,1){34}}
              \put(435,180){\circle{18}}
              \put(396,179){\circle{18}}
              \put(416,199){\circle{18}}
              \put(294,200){\circle{18}}
              \put(313,219){\circle{18}}
              \put(275,219){\circle{18}}
              \put(279,303){\circle{18}}
              \put(279,335){\circle{18}}
              \put(275,420){\circle{18}}
              \put(414,160){\line(1,1){39}}
              \put(373,202){\line(1,-1){41}}
              \put(376,160){\line(1,1){77}}
              \put(254,161){\line(1,1){79}}
              \put(291,238){\line(6,-5){43}}
              \put(255,201){\line(1,1){36}}
              \put(259,319){\line(5,-4){48}}
              \put(308,361){\line(-6,-5){49}}
              \put(296,440){\line(-6,-5){43}}
              \put(254,441){\line(1,-1){40}}
              \put(295,401){\line(0,1){40}}
              \put(374,440){\line(1,0){40}}
              \put(414,439){\line(0,-1){40}}
              \put(414,400){\line(-1,0){40}}
              \put(254,361){\line(1,-1){40}}
              \put(294,319){\line(-1,-1){39}}
              \put(415,321){\line(1,2){20}}
              \put(416,321){\line(1,-2){19}}
              \put(376,361){\line(1,-2){19}}
              \put(395,321){\line(-1,-2){20}}
              \put(254,241){\line(1,-1){79}}
              \put(376,239){\line(1,-1){79}}
              \put(213,203){\line(-1,-1){16}}
              \put(173,163){\line(1,1){16}}
              \put(173,163){\line(-1,1){15}}
              \put(133,203){\line(1,-1){18}}
              \put(135,241){\line(1,-1){79}}
              \put(135,162){\line(1,1){35}}
              \put(179,210){\line(1,1){34}}
              \put(93,203){\line(-1,1){15}}
              \put(52,243){\line(1,-1){17}}
              \put(37,228){\line(1,1){15}}
              \put(13,203){\line(1,1){16}}
              \put(58,209){\line(1,1){34}}
              \put(13,164){\line(1,1){35}}
              \put(13,243){\line(1,-1){79}}
              \put(154,323){\line(-1,-2){20}}
              \put(135,363){\line(1,-2){19}}
              \put(175,323){\line(1,-2){19}}
              \put(174,323){\line(1,2){20}}
              \put(42,301){\line(4,-3){24}}
              \put(17,323){\line(6,-5){16}}
              \put(33,335){\line(-4,-3){16}}
              \put(66,363){\line(-6,-5){24}}
              \put(53,321){\line(-1,-1){39}}
              \put(13,363){\line(1,-1){40}}
              \put(173,402){\line(-1,0){40}}
              \put(173,441){\line(0,-1){40}}
              \put(133,442){\line(1,0){40}}
              \put(28,420){\line(-1,-1){15}}
              \put(54,443){\line(-1,-1){16}}
              \put(54,403){\line(0,1){40}}
              \put(13,443){\line(1,-1){40}}
\end{picture}}

\begin{center} 
{ \bf Figure 1 - Generalized Reidemeister Moves for Virtuals} \end{center}  
\vspace{3mm}

 {\tt    \setlength{\unitlength}{0.92pt}
\begin{picture}(243,311)
\thicklines   \put(161,22){\circle{20}}
              \put(82,21){\circle{20}}
              \put(161,123){\line(0,-1){44}}
              \put(121,124){\line(0,-1){44}}
              \put(80,125){\line(0,-1){44}}
              \put(160,45){\line(0,-1){44}}
              \put(80,45){\line(0,-1){44}}
              \put(202,23){\line(1,1){39}}
              \put(41,22){\line(1,0){161}}
              \put(1,62){\line(1,-1){40}}
              \put(160,226){\line(0,-1){44}}
              \put(80,228){\line(0,-1){44}}
              \put(161,284){\circle{20}}
              \put(122,284){\circle{20}}
              \put(82,284){\circle{20}}
              \put(160,310){\line(0,-1){44}}
              \put(121,310){\line(0,-1){44}}
              \put(81,310){\line(0,-1){44}}
              \put(201,285){\line(1,-1){40}}
              \put(41,285){\line(1,0){160}}
              \put(1,245){\line(1,1){40}}
              \put(119,158){\vector(-1,0){60}}
              \put(103,158){\vector(1,0){76}}
\end{picture}}
 
\begin{center} 
{ \bf Figure 1.1 - Schema for the Detour Move} \end{center}  
\vspace{3mm}

Virtual knot theory is an extension of classical knot theory. In this extension one adds 
a {\em virtual crossing} (See Figure 1) that is neither an over-crossing
nor an under-crossing. We shall refer to the usual diagrammatic crossings, that is those without circles, as {\em real} crossings
to distinguish them from the virtual crossings.
A virtual crossing is represented by two crossing arcs with a small circle
placed around the crossing point. 
\bigbreak

The allowed moves on virtual diagrams are a generalization of the Reidemeister moves for classical knot and link
diagrams. We show the classical Reidemeister moves as part (A) of Figure 1. These classical moves are part of
virtual equivalence where no changes are made to the virtual crossings.
Taken by themselves, the virtual crossings behave as diagrammatic
permutations. Specifically, we have the flat Reidemeister moves (B) for virtual crossings as shown in
Figure 1.  In Figure 1 we also illustrate a basic move (C) that interrelates real and virtual 
crossings. In this move an arc going through a consecutive sequence of two virtual crossings
can be moved across a
single real crossing. In fact, it is consequence of  moves (B) and (C) for virtual crossings that an
arc going through any consecutive
sequence of virtual crossings can be moved anywhere in the diagram keeping the endpoints fixed
and writing the places 
where the moved arc now crosses the diagram as new virtual crossings. This is shown schematically in Figure 1.1.
We call the move in Figure 1.1 the {\em detour}, and note that the detour move is equivalent to 
having all the moves of type (B) and (C) of Figure 1.  This extended move set (Reidemeister moves
plus the detour move or the equivalent moves (B) and (C)) constitutes the move set for virtual knots and links. 
\bigbreak

There is a useful topological interpretation for this virtual theory in terms of embeddings of links
in thickened surfaces. See \cite{VKT,DVK}.  Regard each 
virtual crossing as a shorthand for a detour of one of the arcs in the crossing through a 1-handle
that has been attached to the 2-sphere of the original diagram. The two choices for 
the 1-handle detour are homeomorphic to each other (as abstract surfaces with boundary a circle)
since there is no a priori difference between the meridian and the longitude of a torus.  
By interpreting each virtual crossing in this way, we
obtain an embedding of a collection of circles into a thickened surface  $S_{g} \times R$ where $g$ is the 
number of virtual crossings in the original diagram $L$, $S_{g}$ is a compact oriented surface of genus $g$
and $R$ denotes the real line.  
Thus to each virtual diagram $L$ we 
obtain an embedding $s(L) \longrightarrow S_{g(L)} \times R$ where $g(L)$ is the number of virtual 
crossings of $L$ and $s(L)$ is a disjoint union of circles.  We say that two such surface embeddings are
{\em stably equivalent} if one can be obtained from another by isotopy in the thickened surfaces, 
homeomorphisms of the surfaces and the addition or subtraction of empty handles.  Then we have the
\smallbreak
\noindent
{\bf Theorem \cite{VKT,DKT}.} {\em Two virtual link diagrams are equivalent if and only if their correspondent 
surface embeddings are stably equivalent.}  
\smallbreak
\noindent
In \cite{VKT} this result is sketched. A complete proof will appear in \cite{DKT}.
The surface embedding interpretation of virtuals is useful since it converts their equivalence to a
topological question. However, the stabilization makes classification difficult since
one cannot rely on any single surface embedding. The diagrammatic version 
of virtuals embodies the stabilization in the detour moves. We shall rely on the diagrammatic 
approach here.
\bigbreak

\section{Biquandles and a Generalized Alexander Polynomial $G_{K}(s,t)$}

The {\em biquandle} \cite{DVK, FJK, CS} is an algebra associated with the diagram
that is invariant (up to isomorphism) under the generalized Reidemeister moves for virtual knots and links.
The operations in this algebra are motivated by the formation of labels for the edges of the diagram and the intended 
invariance under the moves. We will give the abstract definition of the biquandle after a discussion of these knot
theoretic issues. View Figure 2.  In this Figure we have shown the format for the operations in a biquandle. The overcrossing arc
has two labels, one on each side of the crossing. In a {\em biquandle}
there is an algebra element labeling  each {\em edge} of the diagram.  An edge of the diagram corresponds to 
an edge of the underlying plane graph of that diagram.
\bigbreak

Let the edges oriented toward a crossing in a diagram be called the {\em input} edges for the crossing, and 
the edges oriented away from the crossing be called the {\em output} edges for the crossing.
Let $a$ and $b$ be the input edges for a positive crossing, with $a$ the label of the undercrossing
input and $b$ the label on the overcrossing input. Then in the biquandle, we label the undercrossing
output by $$c=a^{b}$$ just as in the case of the quandle, but the overcrossing output is labeled
$$d= b_{a}.$$ We usually read $a^{b}$ as -- the undercrossing line $a$ is acted upon by the 
overcrossing line $b$ to produce the output $c=a^{b}.$ In the same way, we can read $b_{a}$ as -- 
the overcossing line $b$ is operated on by the undercrossing line $a$ to produce the output 
$d= b_{a}.$  The biquandle labels for a negative crossing are similar but with 
an overline (denoting an operation of order two) placed on the letters. 
Thus in the case of the negative crossing, we would write

 $$c=a^{\o{b}} \quad \mbox{and} \quad d=b_{\o{a}}.$$
 
\noindent To form the biquandle, $BQ(K)$, we take one generator for each edge of the diagram and two 
relations at each crossing (as described above). This system of generators and relations is then 
regarded as encoding an algebra that is generated freely by the biquandle operations as 
concatenations of these symbols and subject to the biquandle algebra axioms.  These axioms (which we will describe below)
are a transcription 
in the biquandle language of the requirement that this algebra be invariant under Reidemeister moves
on the diagram. 
\vspace{30mm}

{\tt    \setlength{\unitlength}{0.92pt}
\begin{picture}(280,163)
\thicklines   \put(131,162){\line(0,-1){159}}
              \put(1,1){\framebox(278,161){}}
              \put(231,104){\makebox(21,19){$a^{\o{b}} = a\, \UL{b}$}}
              \put(167,59){\makebox(20,15){$b_{\o{a}} = b\, \LL{a}$}}
              \put(233,62){\makebox(20,16){$b$}}
              \put(212,38){\makebox(18,17){$a$}}
              \put(73,106){\makebox(23,23){$a^{b} = a\, \UR{b}$}}
              \put(77,56){\makebox(32,22){$b_{a} = b\, \LR{a}$}}
              \put(10,86){\makebox(20,16){$b$}}
              \put(52,43){\makebox(18,17){$a$}}
              \put(211,89){\vector(0,1){33}}
              \put(211,43){\vector(0,1){32}}
              \put(250,82){\vector(-1,0){78}}
              \put(51,90){\vector(0,1){33}}
              \put(51,43){\vector(0,1){34}}
              \put(11,83){\vector(1,0){80}}
\end{picture}}

\begin{center}
{\bf Figure 2 - Biquandle Relations at a Crossing} \end{center}  \vspace{3mm}

\noindent Another way to write this formalism for the biquandle is as follows

$$a^{b} = a\, \UR{b}$$
$$a_{b} = a\, \LR{b}$$
$$a^{\o{b}} = a\, \UL{b}$$
$$a_{\o{b}} = a\, \LL{b}.$$

\noindent We call this the {\em operator formalism} for the biquandle. The operator formalism
has advantages when one is performing calculations, since it it possible to maintain the formulas
on a line rather than extending them up and down the page as in the exponential notation. On the 
other hand the exponential notation has intuitive familiarity and is good for displaying certain 
results. The axioms for the biquandle, are exactly the rules needed for invariance of this structure
under the Reidemeister moves. Note that in analyzing invariance under Reidemeister moves, we visualize
representative parts of link diagrams with biquandle labels on their edges. The primary labeling
occurs at a crossing. At a positive crossing with over input $b$ and under input $a$, the under output 
is  labeled $a\,\UR{b}$ and the over output is labeled $b\,\LR{a}\,$.   
At a negative crossing with over input $b$ and under input $a$, the under output 
is  labeled $a\,\UL{b}$ and the over output is labeled $b\,\LL{a}$. At a virtual crossing there is
no change in the labeling of the lines that cross one another.  
\bigbreak

\noindent{\bf Remark.} A remark is in order about the relationship of the operator notations with the usual conventions for
binary algebraic operations.  Let $a*b = a^{b} = a\, \UR{b}.$ We are asserting that the biquandle comes equipped with four
binary operations of which one is $a*b.$  Here is how these notations are related to the usual parenthesizations:
\begin{enumerate}
\item $(a*b)*c = (a^{b})^{c} = a^{bc} = a\, \UR{b}\, \UR{c}$
\item $a*(b*c) = a^{b^{c}} = a\, \UR{b\, \UR{c}}$
\end{enumerate}
\noindent From this the reader should see that the exponential and operator notations allow us to express biquandle equation with 
a minimum of parentheses.
\bigbreak

In Figure 3 we illustrate
the effect of these conventions and how it leads to the following algebraic transcription of the 
directly oriented second Reidemeister move:

$$a = a\,\UR{b}\,\UL{b\,\LR{a}} \quad \mbox{or} \quad a = a^{b \o{b_{a}}},$$
$$b = b\,\LR{a}\,\LL{a\,\UR{b}} \quad \mbox{or} \quad b= b_{a \o{a^{b}}}.$$

\vspace{5mm}

{\tt    \setlength{\unitlength}{0.92pt}
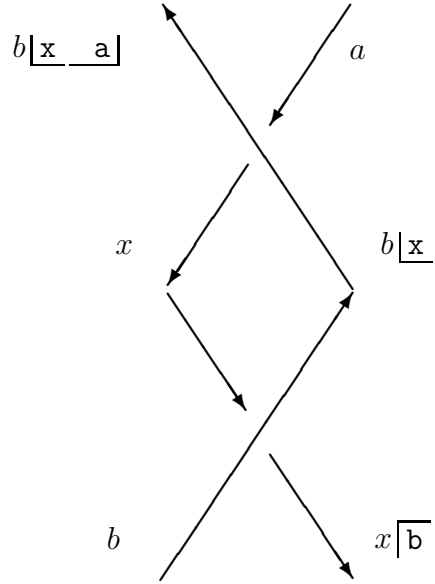
\begin{picture}(323,224)
\thicklines   \put(201,222){\makebox(121,40){$a = a\,\UR{b}\,\UL{b\,\LR{a}}$}}
              \put(41,123){\makebox(81,40){$a\,\UR{b}$}}
              \put(200,122){\makebox(82,41){$b\,\LR{a}$}}
              \put(201,2){\makebox(40,41){$a$}}
              \put(81,2){\makebox(41,41){$b$}}
              \put(1,222){\makebox(121,41){$b = b\,\LR{a}\,\LL{a\,\UR{b}}$}}
              \put(164,186){\vector(2,3){36}}
              \put(123,123){\vector(2,3){34}}
              \put(158,69){\vector(-2,3){36}}
              \put(201,2){\vector(-2,3){36}}
              \put(200,123){\vector(-2,3){78}}
              \put(121,3){\vector(2,3){79}}
\end{picture}}

\begin{center} { \bf Figure 3 --- Direct Two Move } \end{center} 
\vspace{3mm}

{\tt    \setlength{\unitlength}{0.92pt}
\begin{picture}(242,350)
\thicklines   \put(41,1){\makebox(160,42){$\exists x \ni x = a\,\UR{b\,\LL{x}} \quad \mbox{,} \quad a = x\,\UL{b} \quad \mbox{and} \quad b= b\,\LL{x}\, \LR{a} $}}
              \put(1,281){\makebox(80,42){$b\,\LL{x}\, \LR{a}$}}
              \put(26,200){\makebox(78,42){$x$}}
              \put(141,200){\makebox(80,42){$b\,\LL{x}$}}
              \put(139,80){\makebox(81,42){$x\,\UL{b}$}}
              \put(40,82){\makebox(41,39){$b$}}
              \put(141,282){\makebox(40,39){$a$}}
              \put(125,136){\vector(2,-3){34}}
              \put(83,202){\vector(2,-3){32}}
              \put(116,255){\vector(-2,-3){33}}
              \put(158,321){\vector(-2,-3){33}}
              \put(159,204){\vector(-2,3){78}}
              \put(80,84){\vector(2,3){79}}
\end{picture}}

\begin{center} { \bf Figure 4 --- Reverse Two Move } \end{center}
\vspace{3mm}

\noindent The reverse oriented second Reidemeister move gives a different sort of identity, as shown
in Figure 4. For the reverse oriented move, we must assert that given elements $a$ and $b$
in the biquandle, then there exists an element $x$ such that 

$$x = a\,\UR{b\,\LL{x}} \quad \mbox{,} \quad a = x\,\UL{b} \quad \mbox{and} \quad b= b\,\LL{x}\, \LR{a}\,.$$
\bigbreak

By reversing the arrows in Figure 4 we obtain a second statement for invariance under the type two move, saying the same thing with
the operations reversed: Given elements $a$ and $b$
in the biquandle, then there exists an element $x$ such that 

$$x = a\,\UL{b\,\UR{x}} \quad \mbox{,} \quad a = x\,\UR{b} \quad \mbox{and} \quad b= b\,\UR{x}\, \LL{a}\,.$$

\noindent There is no neccessary relation between the $x$ in the first statement and the $x$ in the second statement. 
\bigbreak

\noindent
These assertions about the existence of $x$ can be viewed as asserting the existence of 
fixed points for a certain operators. In this case such an operator is $F(x) =a\,\UR{b\,\LL{x}}\,$.
It is characteristic of certain axioms in the biquandle that they demand the existence of such 
fixed points. Another example is the axiom corresponding to the first Reidemeister move 
(one of them) as illustrated in Figure 5. This axiom states that given an element $a$ in the 
biquandle, then there exists an $x$ in the biquandle such that $x=a \,\LR{x}$ and that 
$a = x \,\UR{a}$. In this case the operator is $G(x) = a\,\LR{x}\,$.

{\tt    \setlength{\unitlength}{0.92pt}
\begin{picture}(323,246)
\thicklines   \put(82,3){\makebox(122,40){$\exists x \ni x=a \,\LR{x} \quad \mbox{and} \quad a = x \,\UR{a}$}}
              \put(224,162){\makebox(78,40){$x \,\UR{a}$}}
              \put(183,82){\makebox(79,40){$a \,\LR{x}$}}
              \put(61,81){\makebox(43,42){x}}
              \put(21,162){\makebox(41,41){a}}
              \put(171,148){\vector(4,3){71}}
              \put(83,82){\vector(4,3){74}}
              \put(122,43){\vector(-1,1){39}}
              \put(201,123){\vector(-1,-1){79}}
              \put(43,203){\vector(2,-1){158}}
\end{picture}}

\begin{center} { \bf Figure 5 --- First Move } \end{center} 
\vspace{3mm}

It is unusual that an algebra
would have axioms asserting the existence of fixed points with respect to operations involving its
own elements. We plan to take up the study of this aspect of biquandles in a separate publication.
\bigbreak

The biquandle relations for invariance under the third Reidemeister move are shown in Figure 6.
The version of the 
third Reidemeister move shown in this figure yields the algebraic relations:

$$a\, \UR{b} \,\UR{c} = a\, \UR{c\, \LR{b}} \,\UR{b\, \UR{c}} \quad \mbox{or} \quad a^{b c} = a^{c_{b} b^{c}},$$

$$c\, \LR{b} \,\LR{a} = c\, \LR{a\, \UR{b}} \,\LR{b\, \LR{a}} \quad \mbox{or} \quad c_{b a} = c_{a^{b} b_{a}},$$

$$b\, \LR{a} \, \UR{c\, \LR{a\, \UR{b}}} = b\, \UR{c}\, \LR{a\, \UR{c\, \LR{b}}} \quad \mbox{or} \quad
 (b_{a})^{c_{a^{b}}} = (b^{c})_{a^{c_{b}}}.$$

{\tt    \setlength{\unitlength}{0.92pt}
\begin{picture}(173,371)
\thicklines   \put(120,63){\vector(1,1){26}}
              \put(87,29){\vector(1,1){25}}
              \put(63,53){\vector(1,-1){24}}
              \put(26,92){\vector(1,-1){28}}
              \put(121,314){\vector(1,-1){25}}
              \put(84,350){\vector(1,-1){28}}
              \put(59,324){\vector(1,1){25}}
              \put(26,290){\vector(1,1){25}}
              \put(147,148){\vector(-1,-1){120}}
              \put(27,148){\vector(1,-1){56}}
              \put(92,82){\vector(1,-1){55}}
              \put(91,283){\vector(1,-1){55}}
              \put(26,349){\vector(1,-1){56}}
              \put(146,349){\vector(-1,-1){120}}
\thinlines    \put(196,2){\makebox(23,23){$b\, \UR{c}\, \LR{a\, \UR{c\, \LR{b}}}$}}
              \put(76,5){\makebox(22,21){$a\, \UR{c\, \LR{b}}$}}
              \put(0,1){\makebox(25,24){$c\, \LR{b} \,\LR{a}$}}
              \put(102,72){\makebox(26,24){$b\, \UR{c}$}}
              \put(48,72){\makebox(23,24){$c\, \LR{b}$}}
              \put(196,77){\makebox(22,24){$a\, \UR{c\, \LR{b}} \,\UR{b\, \UR{c}}$}}
              \put(149,138){\makebox(23,25){$c$}}
              \put(5,139){\makebox(21,23){$b$}}
              \put(4,79){\makebox(20,21){$a$}}
              \put(144,200){\makebox(22,23){$b\, \LR{a} \, \UR{c\, \LR{a\, \UR{b}}}$}}
              \put(6,203){\makebox(25,23){$c\, \LR{a\, \UR{b}} \,\LR{b\, \LR{a}}$}}
              \put(164,279){\makebox(22,25){$a\, \UR{b} \,\UR{c}$}}
              \put(104,280){\makebox(21,23){}}
              \put(48,281){\makebox(23,21){$b\, \LR{a}$}}
              \put(149,340){\makebox(22,21){c}}
              \put(74,349){\makebox(21,21){$a\, \UR{b}$}}
              \put(3,349){\makebox(23,21){$b$}}
              \put(1,282){\makebox(20,20){$a$}}
\end{picture}}

\begin{center} { \bf Figure 6 --- Third Move} \end{center}
\vspace{3mm}

\noindent
The reader will note that if we replace the diagrams of Figure 6 with diagrams with all negative crossings then we will get a second
triple of equations identical to the above equations but with all right operator symbols replaced by the corresponding left operator
symbols (equivalently -- with all exponent literals replaced by their barred versions).  Here are the operator versions of these 
equations. We refrain from writing the exponential versions because of the prolixity of barred variables.

$$a\, \UL{b} \,\UL{c} = a\, \UL{c\, \LL{b}} \,\UL{b\, \UL{c}},$$

$$c\, \LL{b} \,\LL{a} = c\, \LL{a\, \UL{b}} \,\LL{b\, \LL{a}},$$

$$b\, \LL{a} \, \UL{c\, \LL{a\, \UL{b}}} = b\, \UL{c}\, \LL{a\, \UL{c\, \LL{b}}}.$$

\vspace{3mm}

We now have a complete set of axioms, for it is a fact \cite{KNOTS} that the third Reidemeister move with the orientation shown in 
Figure 6 and either all positive crossings (as shown in that Figure) or all negative crossings, is sufficient to generate all the
other cases of third Reidemeister move just so long as we have both oriented forms of the second Reidemeister move.
Consequently, we can now give the full definition of the biquandle.
\bigbreak

\noindent {\bf Definition.}  A {\em biquandle} $B$ is a set with four binary operations indicated by the conventions we have explained
above:  $a^{b} \,\mbox{,} \, a^{\o{b}} \, \mbox{,} \,  a_{b} \,\mbox{,} \, a_{\o{b}}.$ We shall refer to the operations with barred 
variables as the {\em left} operations and the operations without barred variables as the {\em right} operations. The biquandle is 
closed under these operations and the following axioms are satisfied:

\begin{enumerate}
\item   For any elements $a$ and $b$ in $B$ we have 

$$a = a^{b \o{b_{a}}}  \quad \mbox{and} \quad  b= b_{a \o{a^{b}}} \quad \mbox{and}$$

$$a = a^{\o{b}b_{\o{a}}}  \quad \mbox{and} \quad  b= b_{\o{a} a^{\o{b}}}.$$

\item  Given elements $a$ and $b$
in $B$, then there exists an element $x$ such that 

$$x = a^{b_{\o{x}}} \mbox{,} \quad a = x^{\o{b}} \quad \mbox{and} \quad b= b_{\o{x}a}.$$

\noindent Given elements $a$ and $b$
in $B$, then there exists an element $x$ such that 

$$x = a^{\o{b_{x}}} \mbox{,} \quad a = x^{b} \quad \mbox{and} \quad b= b_{x\o{a}}.$$

\item For any $a$ , $b$ , $c$ in $B$ the following equations hold and the same equations hold when all right operations are 
replaced in these equations by left operations.

$$a^{b c} = a^{c_{b} b^{c}} \mbox{,} \quad c_{b a} = c_{a^{b} b_{a}} \mbox{,} \quad (b_{a})^{c_{a^{b}}} = (b^{c})_{a^{c_{b}}}.$$

\item Given an element $a$ in $B$, then there exists an $x$ in the biquandle such that $x=a_{x}$ and  
$a = x^{a}.$ Given an element $a$ in $B$, then there exists an $x$ in the biquandle such that $x=a^{\o{x}}$ and  
$a = x_{\o{a}}.$
\end{enumerate}

These axioms are transcriptions of the Reidemeister moves. The first axiom transcribes the directly oriented second Reidemeister move. 
The second axiom transcribes the reverse oriented Reidemeister move. The third axiom transcribes the third Reidemeister move as we have
described it in Figure 6. The fourth axiom transcribes the first Reidemeister move. Much more work is needed in exploring these 
algebras and their applications to knot theory.
\bigbreak

\subsection{The Alexander Biquandle}

\noindent In order to realize a specific example of a biquandle structure, suppose that 

$$a\,\UR{b} = ta + vb$$
$$a\,\LR{b} =sa + ub$$

\noindent where  $a$,$b$,$c$ are elements of a module $M$ over a ring 
$R$ and $t$,$s$,$v$
and $u$ are in $R$. We use invariance under the Reidemeister moves to 
determine relations among these coefficients. 
\vspace{3mm}

Taking the equation for the third Reidemeister move discussed  above, we have

$$a \UR{b} \,\UR{c} = t(ta+vb) + vc = t^{2}a + tvb +vc$$
$$a \UR{c \LR{b}} \,\UR{b \UR{c}} = t(ta + v(sc + ub)) + v(tb + vc)$$
$$= t^{2}a + tv(u+1)b + v(ts+v)c.$$

\noindent From this we see that we have a solution to the equation for the third Reidemeister move
if $u=0$ and $v=1-st$. Assuming that $t$ and $s$ are invertible, it is not hard to see that 
the following equations not only solve this single Reideimeister move, but they give a biquandle 
structure, satisfying all the moves.

$$a\,\UR{b} = ta + (1-st)b \, \mbox{,} \quad a\,\LR{b} = sa$$
$$a\,\UL{b} = t^{-1}a + (1-s^{-1}t^{-1})b \, \mbox{,} \quad a\,\LL{b} = s^{-1}a.$$

\noindent Thus we have a simple generalization of the Alexander quandle and we shall refer to this 
structure, with the equations given above, as the {\em Alexander Biquandle}. 
\vspace{3mm}

Just as one can define the Alexander Module of a classical knot, we have the Alexander Biquandle of
a virtual knot or link, obtained by taking one generator for each {\em edge} of the knot diagram and
taking the module relations in the above linear form. Let $ABQ(K)$ denote this module structure for an
oriented link $K$. That is, $ABQ(K)$ is the module generated by the edges of the diagram, modulo the submodule
generated by the relations. This module then has a biquandle structure specified by the operations defined above for an
Alexnder Biquandle.  We first construct the module and then note that it has a biquandle structure.
See Figures 7,8 and 9 for an illustration of the Alexander Biquandle labelings at a crossing.
\vspace{3mm}

{\tt    \setlength{\unitlength}{0.92pt}
\begin{picture}(309,475)
\thinlines    \put(217,25){\vector(-1,1){158}}
              \put(58,24){\vector(1,1){71}}
              \put(146,114){\vector(1,1){70}}
              \put(83,1){\framebox(40,41){$a$}}
              \put(227,7){\framebox(41,40){$b$}}
              \put(22,188){\framebox(61,41){$s^{-1}b$}}
              \put(107,188){\framebox(201,41){$t^{-1}a+(1-s^{-1}t^{-1})b$}}
              \put(1,432){\framebox(120,42){$ta+(1-st)b$}}
              \put(222,414){\framebox(39,41){$sb$}}
              \put(16,242){\framebox(40,42){$b$}}
              \put(226,248){\framebox(40,42){$a$}}
              \put(129,354){\vector(-1,1){71}}
              \put(218,265){\vector(-1,1){72}}
              \put(57,265){\vector(1,1){161}}
\end{picture}}
\bigbreak

\begin{center} {\bf Figure 7 - Alexander Biquandle Labeling at a Crossing} \end{center}
\vspace{3mm}

{\tt    \setlength{\unitlength}{0.92pt}
\begin{picture}(344,352)
\thicklines   \put(286,48){\framebox(41,43){$d$}}
              \put(156,206){\framebox(41,42){$c$}}
              \put(209,102){\framebox(42,42){$b$}}
              \put(1,216){\framebox(40,41){$a$}}
              \put(154,91){\circle{30}}
              \put(271,198){\vector(1,-1){71}}
              \put(182,286){\vector(1,-1){75}}
              \put(80,203){\vector(-2,3){55}}
              \put(185,45){\vector(-2,3){86}}
              \put(63,2){\vector(-1,3){42}}
              \put(323,263){\vector(-1,-1){260}}
              \put(212,348){\vector(4,-3){113}}
              \put(25,287){\vector(3,1){186}}
              \put(342,125){\vector(-2,-1){157}}
              \put(22,125){\vector(1,1){160}}
\end{picture}}

\begin{center} {\bf Figure 8 - A Virtual Knot Fully Labeled} \end{center}
\vspace{3mm}

{\tt    \setlength{\unitlength}{0.92pt}
\begin{picture}(333,352)
\thicklines   \put(132,205){\framebox(67,41){$s^{2}a$}}
              \put(204,114){\framebox(41,41){$sa$}}
              \put(291,165){\framebox(41,43){$d$}}
              \put(1,192){\framebox(40,41){$a$}}
              \put(138,91){\circle{30}}
              \put(255,198){\vector(1,-1){71}}
              \put(166,286){\vector(1,-1){75}}
              \put(64,203){\vector(-2,3){55}}
              \put(169,45){\vector(-2,3){86}}
              \put(47,2){\vector(-1,3){42}}
              \put(307,263){\vector(-1,-1){260}}
              \put(196,348){\vector(4,-3){113}}
              \put(9,287){\vector(3,1){186}}
              \put(326,125){\vector(-2,-1){157}}
              \put(6,125){\vector(1,1){160}}
\end{picture}}

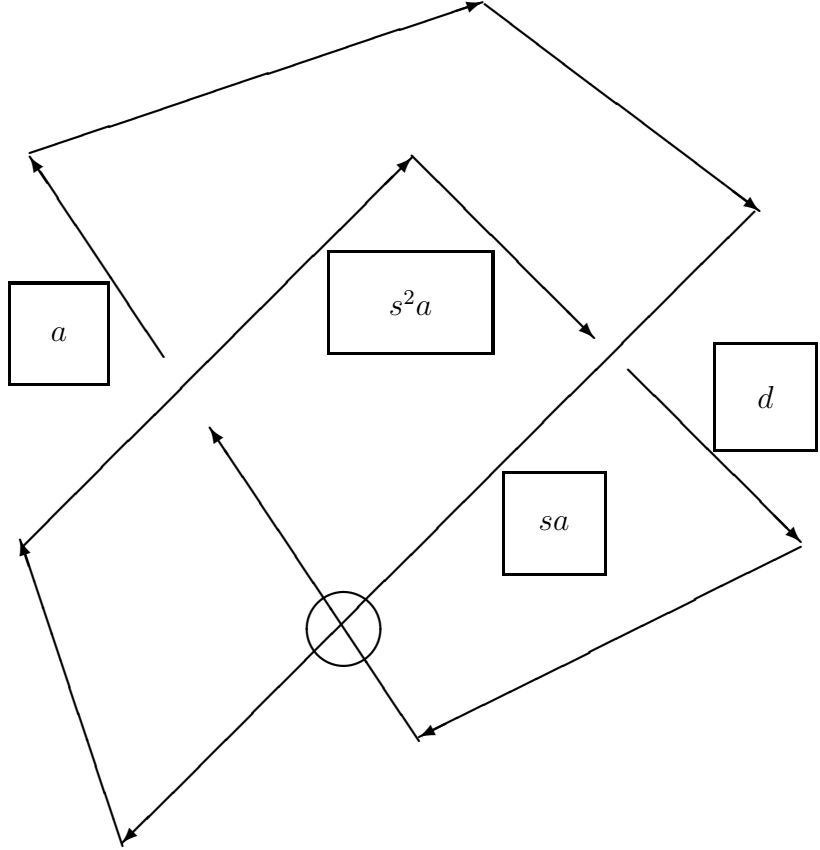
\begin{center} {\bf Figure 9 - A Virtual Knot with Lower Operations Labeled} \end{center}
\vspace{3mm}

For example, consider the virtual knot in Figure 8.  This knot gives rise to a
biquandle with generators $a$,$b$,$c$,$d$ and relations

$$a= d\, \UR{b} \, \mbox{,} \quad c= b\, \LR{d} \, \mbox{,} \quad d= c\, \UR{a} \, \mbox{,} \quad b= a\, \LR{c}.$$

\noindent writing these out in $ABQ(K)$, we have

$$a = td + (1-st)b \, \mbox{,} \quad c=sb \, \mbox{,} \quad d=tc+(1-st)a \, \mbox{,} \quad b=sa.$$

\noindent eliminating $c$ and $b$ and rewriting, we find

$$a=td+(1-st)sa$$
$$d=ts^{2}a+(1-st)a$$

Note that these relations can be written directly from the diagram as indicated in Figure 9 if we perform the lower biquandle operations
directly on the diagram.  This is the most convenient algorithm for producing the relations. 
\smallbreak

\noindent We can write these as a list of relations
 
$$(s-s^{2}t-1)a +td = 0$$
$$(s^{2}t+1-st)a -d = 0$$

\noindent for the Alexander Biquandle as a module over $Z[s,s^{-1}, t, t^{-1}].$
The relations can be expressed concisely with the matrix of coefficients of this system of equations:

 $$M =  \left[
\begin{array}{cc}
     s-s^{2}t-1 & t  \\
     (s^{2}t+1-st) & -1
\end{array}
\right]. $$

\noindent The determinant of $M$ is, up to multiples of $\pm s^{i}t^{j}$ for integers $i$ and $j$, an invariant of the virtual knot 
or link $K$. We shall denote this determinant by $G_{K}(s,t)$ and call it the {\em generalized Alexander polynomial for $K$}.
A key fact about $G_{K}(s,t)$ is that {\em $G_{K}(s,t)=0$ if $K$ is equivalent to a classical diagram}. 
This is seen by noting that 
in a classical diagram one of the relations will be a consequence of the others. 
\smallbreak

\noindent In this case we have

$$G_{K} = (1-s) +(s^{2}-1)t + (s-s^{2})t^{2},$$
\noindent which shows that the knot in question is non-trivial and non-classical.
\bigbreak

Here is another example of the use of this polynomial.  Let $D$ denote the diagram in Figure 10.
It is not hard to see that this virtual knot has unit Jones polynomial, and that the fundamental
group is isomorphic to the integers.  The biquandle does detect the knottedness of $D$.
The relations are

$$ a\, \UR{d} = b, \, d\, \LR{a} = e \, \mbox{,} \quad c\, \UR{e} = d, \, e\,\LR{c} = f \, \mbox{,} \quad f\, \UL{b} = a, \, b\, \LL{f} = c$$

\noindent from which we obtain the relations (eliminating $c$, $e$ and $f$)

$$b = ta +(1-tv)d \, \mbox{,} \quad d = ts^{-1}b + (1-ts)sd \, \mbox{,} \quad a = t^{-1}s^{2}d + (1-t^{-1}s^{-1})b \, .$$

\noindent The determinant of this system is the generalized Alexander polynomial for $D$:

$$t^{2}(s^{2}-1) + t(s^{-1}+1-s-s^{2}) + (s-s^{2}).$$

\noindent This proves that $D$ is a non-trivial virtual knot.
\vspace{5mm}

{\tt    \setlength{\unitlength}{0.92pt}
\begin{picture}(326,328)
\thinlines    \put(292,137){\makebox(33,41){$D$}}
\thicklines   \put(25,242){\vector(0,-1){161}}
              \put(26,85){\vector(3,2){241}}
              \put(65,167){\vector(1,2){79}}
              \put(145,326){\vector(1,-2){121}}
              \put(145,5){\vector(-1,2){54}}
              \put(80,135){\vector(-1,2){15}}
              \put(264,84){\vector(-3,2){111}}
              \put(135,169){\vector(-3,2){109}}
              \put(267,245){\vector(-1,-2){36}}
              \put(222,154){\vector(-1,-2){76}}
              \put(372,291){\circle{0}}
              \put(84,202){\circle{22}}
              \put(207,203){\circle{22}}
              \put(206,121){\circle{22}}
              \put(117,4){\makebox(22,24){$a$}}
              \put(73,144){\makebox(19,21){$b$}}
              \put(180,248){\makebox(22,25){$b$}}
              \put(249,106){\makebox(21,27){$c$}}
              \put(160,123){\makebox(16,25){$c$}}
              \put(111,181){\makebox(19,22){$d$}}
              \put(1,90){\makebox(21,24){$d$}}
              \put(116,124){\makebox(17,20){$e$}}
              \put(179,164){\makebox(14,17){$f$}}
              \put(244,177){\makebox(21,22){$f$}}
\end{picture}}

\begin{center} {\bf  Figure 10 -- Unit Jones, Integer Fundamental Group } \end{center}
\vspace{3mm}

In fact the polynomial that we have computed is the same as the polynomial invariant of virtuals 
of Sawollek \cite{SAW} and defined by an alternative method by Silver and Williams \cite{SW}. 
Sawollek defines a module structure essentially the same as our Alexander Biquandle. 
Silver and Williams first define a group.
The Alexander Biquandle proceeds from taking the abelianization of the Silver-Williams group.
\vspace{3mm}

We end this discussion of the Alexander Biquandle with two examples that show clearly its 
limitations. View Figure 11. In this Figure we illustrate two diagrams labeled $K$ and $KI.$
It is not hard to calculate that both $G_{K}(s,t)$ and $G_{KI}(s,t)$ are equal to zero. However,
The Alexander Biquandle of $K$ is non-trivial -- calculation shows that it is isomorphic to the 
free module over $Z[s, s^{-1},t, t^{-1}]$ generated by elements $a$ and $b$ subject to the relation
$(s^{-1} - t -1)(a-b) =0.$ Thus $K$ represents a non-trivial virtual knot. This shows that it is 
possible for a non-trivial virtual diagram to be a connected sum of two trivial virtual diagrams,
and it shows that the Alexander Biquandle can sometimes be more powerful than the polynomial invariant
$G.$ However, the diagram $KI$ also has trivial Alexander Biquandle.  In fact $KI$, discovered by
Kishino \cite{P} is not yet proved to be a knotted virtual.  We conjecture that $KI$ is knotted and 
that its general Biquandle is non-trivial.
\bigbreak

{\tt    \setlength{\unitlength}{0.92pt}
\begin{picture}(272,368)
\thicklines   \put(243,64){\makebox(25,21){$KI$}}
              \put(244,260){\makebox(27,25){$K$}}
              \put(165,45){\line(-1,0){14}}
              \put(82,45){\line(1,0){58}}
              \put(163,123){\line(-1,0){14}}
              \put(82,123){\line(1,0){59}}
              \put(144,163){\line(0,-1){159}}
              \put(3,124){\line(1,-1){78}}
              \put(3,44){\line(1,1){78}}
              \put(163,124){\line(1,-1){79}}
              \put(163,44){\line(1,1){80}}
              \put(105,164){\line(0,-1){38}}
              \put(105,118){\line(0,-1){68}}
              \put(106,37){\line(0,-1){33}}
              \put(145,163){\line(5,-2){99}}
              \put(44,83){\circle{18}}
              \put(204,84){\circle{18}}
              \put(3,123){\line(5,2){103}}
              \put(2,45){\line(5,-2){105}}
              \put(146,5){\line(5,2){97}}
              \put(146,206){\line(5,2){97}}
              \put(2,246){\line(5,-2){105}}
              \put(3,324){\line(5,2){103}}
              \put(204,285){\circle{18}}
              \put(44,284){\circle{18}}
              \put(145,364){\line(5,-2){99}}
              \put(144,238){\line(0,-1){33}}
              \put(144,317){\line(0,-1){69}}
              \put(143,365){\line(0,-1){37}}
              \put(106,238){\line(0,-1){33}}
              \put(105,316){\line(0,-1){68}}
              \put(105,365){\line(0,-1){38}}
              \put(82,245){\line(1,0){81}}
              \put(82,324){\line(1,0){81}}
              \put(163,245){\line(1,1){80}}
              \put(163,325){\line(1,-1){79}}
              \put(3,245){\line(1,1){78}}
              \put(3,325){\line(1,-1){78}}
\end{picture}}

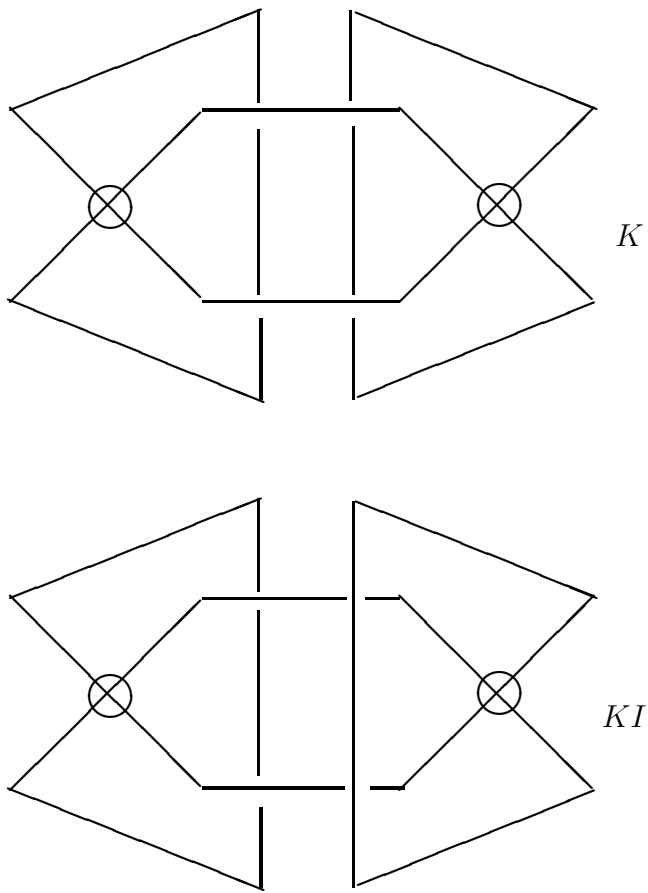
\begin{center} {\bf  Figure 11 -- The Knot $K$ and the Kishino Diagram $KI$ } \end{center}
\vspace{3mm}

\section{A Quantum Model for $G_{K}(s,t)$}

It is our intent in this paper to analyse the structure of the invariant $G_{K}(s,t)$ by rewriting it as a quantum 
invariant and then analysing its state summation. We shall show how the quantum invariant and state sum fit 
into the context of oriented quantum 
algebras. The quantum model for this invariant is obtained in a fashion analogous to the construction of a 
quantum model of the Alexander polynomial in \cite{KS} and  \cite{JKS}. The strategy in those papers was to take the basic 
two dimensional matrix of the Burau representation, view it as a linear transformation  $T: V \longrightarrow V$ on 
a two dimensional module $V$, and them take the induced linear transformation 
$\hat{T}: \Lambda^{*}V \longrightarrow  \Lambda^{*}V$ on the exterior algebra of $V$. This gives a transformation on a
four dimensional module that is a solution to the Yang-Baxter equation.  This solution of the Yang-Baxter equation
then becomes the building block for the corresponding quantum invariant.  In the present instance, we have a generalization
of the Burau representation, and this same procedure can be applied to it.
\bigbreak

The generalized Burau matrix is given by the formula

 $$B =  \left[
\begin{array}{cc}
     1-st & s  \\
     t & 0
\end{array}
\right]$$

\noindent with the inverse matrix

 $$B^{-1} =  \left[
\begin{array}{cc}
     0 & t^{-1}  \\
     s^{-1} & 1-s^{-1}t^{-1}
\end{array}
\right]. $$

\noindent The formulas for $B$ and $B^{-1}$ are easily seen by reference to Figure 7.
\smallbreak

We may regard $B$ as acting on a module $V$ with basis $\{ e_{1}, e_{2} \}$ over $Z[s,s^{-1},t,t^{-1}]$  via the equations

$$Be_{1} = (1-st)e_{1} + te_{2},$$ 
$$Be_{2} = se_{1}.$$

\noindent Letting $\hat{B}$ denote the extension of $B$ to the exterior algebra on $V$, we have
 
$$\hat{B}1 = 1,$$
$$Be_{1} = (1-st)e_{1} + te_{2},$$ 
$$Be_{2} = se_{1},$$
$$\hat{B}e_{1} \wedge e_{2} = Det(B) e_{1} \wedge e_{2} = -st e_{1} \wedge e_{2}.$$ 

\noindent Let $R$ denote the matrix of $\hat{B}$ with respect to the basis $\{1, e_{1}, e_{2}, e_{1} \wedge e_{2} \}$ of 
the exterior algebra $\Lambda^{*}V.$ The matrix $R$ is the matrix of the transformation on the exterior algebra that is
induced from the generalized Burau matrix.

$$R =  \left[
\begin{array}{cccc}
     1 & 0 & 0 & 0 \\
     0 & 1-st & s & 0\\
     0 & t & 0 & 0 \\
     0 & 0 & 0 & -st
\end{array}
\right]$$

\noindent
$R$ is a  $4 \times 4$ matrix solution to the Yang-Baxter equation.  Its inverse $\bar{R}$ is shown below.

$$\bar{R} =  \left[
\begin{array}{cccc}
     1 & 0 & 0 & 0 \\
     0 & 0 & t^{-1} & 0\\
     0 & s^{-1} & 1-s^{-1}t^{-1} & 0 \\
     0 & 0 & 0 & -s^{-1}t^{-1}
\end{array}
\right]$$

\noindent
In our case, we also need the induced transformation for the virtual crossing.  At the level of the generalized Burau 
representation, the matrix for the virtual crossing is the $2 \times 2$ matrix for a transposition:

$$\eta =  \left[
\begin{array}{cc}
     0 & 1  \\
     1 & 0
\end{array}
\right]. $$

\noindent
The corresponding matrix induced on the exterior algebra is

$$\hat{\eta} =  \left[
\begin{array}{cccc}
     1 & 0 & 0 & 0 \\
     0 & 0 & 1 & 0\\
     0 & 1 & 0 & 0 \\
     0 & 0 & 0 & -1
\end{array}
\right]. $$

\noindent This matrix $\hat{\eta}$ is a solution to the Yang-Baxter equation whose square is equal to the identity.
It is this operator that will correspond to the virtual crossings in our quantum invariant for virtuals.
\smallbreak

Now it is convenient to make some changes of variables for this model. We replace $s$ by $\sigma^{2}$ and $t$ by
$\tau^{-2}$. We then replace $R$ by $\sigma^{-1}\tau R$ and $\bar{R}$ by $\sigma \tau^{-1} \bar{R}.$ The result is
the new matrices, shown below where: 

$$z = \sigma^{-1} \tau - \sigma \tau^{-1}$$

 $$R =  \left[
\begin{array}{cccc}
     \sigma^{-1} \tau & 0 & 0 & 0 \\
     0 & z & \sigma \tau & 0\\
     0 & \sigma^{-1} \tau^{-1} & 0 & 0 \\
     0 & 0 & 0 & -\sigma \tau^{-1}
\end{array}
\right]. $$

 $$\bar{R} =  \left[
\begin{array}{cccc}
     \sigma \tau^{-1} & 0 & 0 & 0 \\
     0 & 0 & \sigma \tau & 0\\
     0 & \sigma^{-1} \tau^{-1} & -z & 0 \\
     0 & 0 & 0 & -\sigma^{-1} \tau
\end{array}
\right]. $$

These matrices plus a choice of cup and cap matrices will define the matrix model for this quantum invariant. See \cite{KRCAT} for 
a description of the matrix models. See \cite{KNOTS} for a related discussion of a state summation for the Alexander polynomial.
Here the cup and cap matrices are given by the formula:

$$M_{ab} = (\sqrt{i})^{a}\delta_{ab}$$

\noindent where $i^{2} = -1$ and the matrix $M$ is associated with the clockwise-turning cap and with the clockwise-turning cup, while
the matrix $M^{-1}$ is associated with the counterclockwise-turning cap and with the counterclockwise-turning cup. Here the matrix
indices are $-1$ and $+1$ so that an isolated clockwise loop evaluates as $i^{+1} + i^{-1} = 0$. (Remember that this invariant 
vanishes on classical links.) It is easy to verify that this model is invariant under all but the first Reidemeister moves (classical and flat).
Let $W(K)$ denote this evaluation.  It is then not hard to see that the following normalization $Z(K)$ creates a function on virtual knots 
that is invariant under all of the (generalized) Reidemeister moves:

$$Z(K) = (\sigma^{-1} \tau i)^{rot(K) - v(K)} i^{v(K)} W(K)$$

\noindent where $rot(K)$ denotes the sum of the  Whitney degrees of the underlying plane curves of $K$ (where $rot(K) = 1$ when $K$ is a simple 
clockwise circle in the plane), and $v(K)$ denotes the number of virtual crossings in the diagram $K$.

We can formulate this invariant as a state summation by using the formulas in Figure 12 to expand a given diagram into a sum of evaluations of labeled
signed loops in the plane. Each loop has only virtual crossings and the rule for expanding the virtual tells us that only the $++$ signing receives a minus 
one as vertex weight. Note the distinction between the vertex weights (the matrix entries) and the signs (the matrix indices).
In all other cases the virtual crossing contributes one as a vertex weight. Once there is such a labeled collection of signed loops, it can be tested 
for {\em validity}. A collection of loops is {\em valid} if there is no contradiction in the signs (each curve is uniquely labeled plus or minus). An invalid
labeling is a zero summand in the state sum.  Then each valid labeled signed collection is evaluated by taking the product of the vertex weights multiplied 
by the product of the evaluations of the signed loops. The evaluation of a given signed loop $\lambda$ is equal to $i^{\epsilon(\lambda) rot(\lambda)}$ 
where $\epsilon(\lambda)$ is the sign of the loop and $rot(\lambda)$ is the Whitney degree of the loop.  Note that because of the presence of the
virtual crossings, the Whitney degree of a loop can be any integer.  The state summation evaluates the unnormalized invariant $W(K).$
\smallbreak

\noindent Both $W(K)$ and $Z(K)$ satisfy a skein relation that is just like that of the Conway polynomial:

$$W(K_{+}) - W(K_{-}) = zW(K_{0})$$

$$Z(K_{+}) - Z(K_{-}) = zZ(K_{0})$$

\noindent where $K_{+}$, $K_{-}$ and $K_{0}$ denote three diagrams that differ at one site, with a positive crossing, a negative crossing and a smoothing
respectively. This skein relation, illustrated in Figure 12,
is useful in relative computations, but there are infinitely many virtual links whose evaluation cannot be decided by the
skein relation alone. The simplest example is the ``virtual Hopf link" $H$ of Figure 13. In $H$ there is one crossing and one virtual crossing. Switching this
crossing does not simplify the link. In Figure 13 we illustrate the state sum calculation of $W(H)$. The Figure shows the four contributing states.
It is clear from the Figure that 
$$W(H) = \tau\sigma^{-1} + \tau^{-1}\sigma -\sigma\tau - \sigma^{-1}\tau^{-1} = (\tau - \tau^{-1})(\sigma^{-1} - \sigma)$$
\smallbreak

\noindent This calculation is obtained by evaluating the signed loops according to the description given above.  In particular, the case of the single loop
with plus sign has a rotation number of zero and a multiplicative vertex weight of minus one contributed from the virtual crossing. In all the other cases
the vertex weight of the virtual crossing is plus one. The first two states have rotation number zero, while the second two contribute $i^{-2}=-1$ and
$i^{2} = -1$ respectively. This explains the signs in the formula for $W(H)$.  To obtain $Z(H)$, we note that $rot(H)=0$ and that $H$ has a single virtual 
crossing.  So

$$Z(H) = (\sigma^{-1} \tau i)^{rot(H) - v(H)} i^{v(H)} W(H) = (\sigma^{-1} \tau i)^{0 - 1} i^{1} W(H)$$

$$= \sigma \tau^{-1} W(H) = \sigma \tau^{-1} (\tau - \tau^{-1})(\sigma^{-1} - \sigma) = (1-\tau^{-2})(1 - \sigma^{2}) = (1-t)(1-s).$$

\noindent It is easy to see that this corresponds to the calculation of $G_{H}(s,t)$ by the determinant function.  
The state summation $Z(H)$ provides a normalized
value for the invariant that can be used for skein calculations.
\bigbreak

 The basic result behind the correspondence of $G_{K}(s,t)$ and $Z(K)$ is the 
 \bigbreak
 
 \noindent {\bf Theorem.} {\em For a (virtual) link $K$, the invariants $Z(K)(\sigma = \sqrt{s}, \tau = 1/\sqrt{t})$ and $G_{K}(s,t)$ are equal up to 
 a multiple of $\pm s^{n}t^{m}$ for integers $n$ and $m$ (this being the well-definedness criterion for $G$).}
 \smallbreak
 
 We omit the proof of this result. The argument is the same as that for the relation between the classical Alexander polynomial and its quantum 
 counterpart. See \cite{KS} and \cite{KNOTS}.
 \vspace{80mm}

{\tt    \setlength{\unitlength}{0.92pt}
\begin{picture}(373,106)
\thinlines    \put(174,3){\makebox(80,39){$+\sigma^{-1} \tau^{-1}$}}
              \put(75,2){\makebox(56,41){$+\sigma \tau$}}
              \put(258,60){\makebox(71,40){$-\tau^{-1} \sigma$}}
              \put(139,60){\makebox(70,42){$+\tau \sigma^{-1}$}}
              \put(46,71){\makebox(42,25){$= z$}}
              \put(262,39){\framebox(29,17){$+-$}}
              \put(140,41){\framebox(29,16){$-+$}}
              \put(336,74){\framebox(29,18){$++$}}
              \put(220,75){\framebox(30,17){$--$}}
              \put(100,72){\framebox(31,22){$-+$}}
\thicklines   \put(256,2){\vector(1,1){39}}
              \put(295,2){\vector(-1,1){39}}
              \put(134,5){\vector(1,1){39}}
              \put(173,5){\vector(-1,1){39}}
              \put(329,65){\vector(0,1){39}}
              \put(370,65){\vector(0,1){39}}
              \put(255,64){\vector(0,1){39}}
              \put(214,64){\vector(0,1){39}}
              \put(135,63){\vector(0,1){42}}
              \put(95,64){\vector(0,1){39}}
              \put(18,86){\vector(-1,1){17}}
              \put(42,63){\vector(-1,1){15}}
              \put(2,63){\vector(1,1){39}}
\end{picture}}
\bigbreak

{\tt    \setlength{\unitlength}{0.92pt}
\begin{picture}(374,108)
\thicklines   \put(25,89){\vector(1,1){17}}
              \put(1,67){\vector(1,1){16}}
              \put(41,67){\vector(-1,1){37}}
\thinlines    \put(175,3){\makebox(80,39){$+\sigma^{-1} \tau^{-1}$}}
              \put(76,2){\makebox(56,41){$+\sigma \tau$}}
              \put(259,60){\makebox(71,40){$-\tau \sigma^{-1}$}}
              \put(140,60){\makebox(70,42){$+\tau^{-1} \sigma$}}
              \put(47,71){\makebox(42,25){$= -z$}}
              \put(263,39){\framebox(29,17){$-+$}}
              \put(141,41){\framebox(29,16){$+-$}}
              \put(337,74){\framebox(29,18){$++$}}
              \put(221,75){\framebox(30,17){$--$}}
              \put(101,72){\framebox(31,22){$+-$}}
\thicklines   \put(257,2){\vector(1,1){39}}
              \put(296,2){\vector(-1,1){39}}
              \put(135,5){\vector(1,1){39}}
              \put(174,5){\vector(-1,1){39}}
              \put(330,65){\vector(0,1){39}}
              \put(371,65){\vector(0,1){39}}
              \put(256,64){\vector(0,1){39}}
              \put(215,64){\vector(0,1){39}}
              \put(136,63){\vector(0,1){42}}
              \put(96,64){\vector(0,1){39}}
\end{picture}}
\bigbreak

{\tt    \setlength{\unitlength}{0.92pt}
\begin{picture}(203,43)
\thicklines   \put(125,7){\makebox(29,27){= z}}
              \put(46,8){\makebox(29,25){-}}
              \put(200,2){\vector(0,1){40}}
              \put(160,3){\vector(0,1){39}}
              \put(103,25){\vector(1,1){16}}
              \put(80,3){\vector(1,1){14}}
              \put(119,2){\vector(-1,1){39}}
              \put(16,26){\vector(-1,1){15}}
              \put(40,2){\vector(-1,1){15}}
              \put(1,2){\vector(1,1){39}}
\end{picture}}
\bigbreak

{\tt    \setlength{\unitlength}{0.92pt}
\begin{picture}(364,58)
\thicklines   \put(202,3){\vector(-1,1){39}}
              \put(163,3){\vector(1,1){39}}
              \put(122,2){\vector(-1,1){40}}
              \put(82,2){\vector(1,1){40}}
\thinlines    \put(289,8){\makebox(28,27){+}}
              \put(210,7){\makebox(27,28){+}}
              \put(130,9){\makebox(28,26){-}}
              \put(47,10){\makebox(30,26){=}}
\thicklines   \put(283,2){\vector(-1,1){39}}
              \put(244,2){\vector(1,1){39}}
              \put(362,3){\vector(-1,1){39}}
              \put(323,3){\vector(1,1){39}}
\thinlines    \put(89,39){\framebox(30,17){$--$}}
              \put(169,38){\framebox(29,18){$++$}}
              \put(250,38){\framebox(29,16){$-+$}}
              \put(329,40){\framebox(29,17){$+-$}}
\thicklines   \put(22,24){\circle{22}}
              \put(43,3){\vector(-1,1){38}}
              \put(1,4){\vector(1,1){38}}
\end{picture}}
\bigbreak

\begin{center}
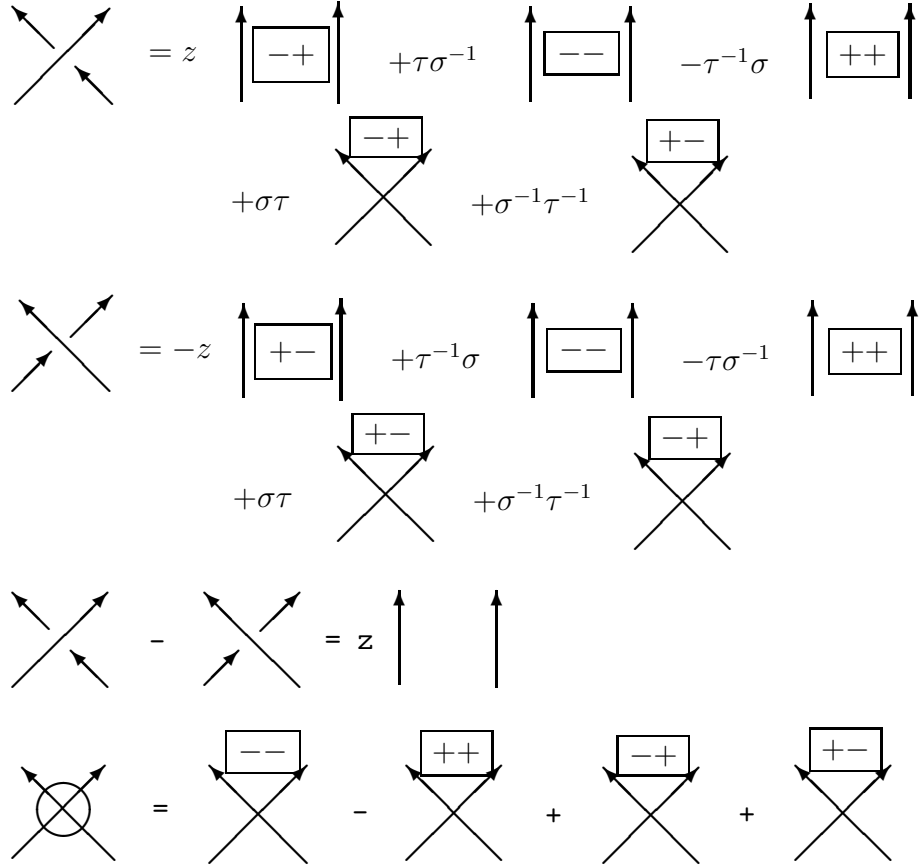
{\bf Figure 12  - Expansion Formulas for the State Summation}\end{center}
\vspace{3mm}

 \noindent{\bf Remark.} It should be remarked that there is a natural multivariable version of the polynomial invariant
$G_{K}(s,t)$ to a polynomial invariant $G_{K}(s, t_{1}, ... , t_{\mu})$ where $\mu$ denotes the number of components in the 
link $K.$ Each link component receives a separate variable, and the biquandle relations at a crossing are determined by the 
label for the undercrossing segment at that crossing. The procedure of extending through the exterior algebra still goes through
to produce a quantum model for this many-variable polynomial.  We will study this generalized invariant in a separate paper.
 \vspace{3mm}

{\tt    \setlength{\unitlength}{0.92pt}
\begin{picture}(402,549)
\thinlines    \put(155,69){\makebox(121,39){$+ \sigma^{-1}\tau^{-1}$}}
              \put(126,187){\makebox(120,41){$+ \sigma\tau$}}
              \put(145,348){\makebox(123,40){$- \tau^{-1}\sigma$}}
              \put(144,466){\makebox(122,42){$= \tau\sigma^{-1}$}}
\thicklines   \put(327,169){\makebox(17,17){+}}
              \put(303,196){\makebox(15,17){-}}
              \put(329,25){\makebox(20,19){-}}
              \put(298,52){\makebox(20,19){+}}
              \put(294,327){\makebox(23,20){+}}
              \put(324,295){\makebox(23,22){+}}
              \put(323,435){\makebox(22,22){-}}
              \put(295,465){\makebox(22,22){-}}
              \put(360,82){\circle{22}}
              \put(398,3){\vector(-1,0){80}}
              \put(398,83){\vector(0,-1){80}}
              \put(318,83){\vector(1,0){80}}
              \put(359,44){\vector(0,1){79}}
              \put(277,44){\vector(1,0){82}}
              \put(278,123){\vector(0,-1){80}}
              \put(358,123){\vector(-1,0){80}}
              \put(316,3){\vector(0,1){80}}
              \put(360,502){\circle{22}}
              \put(398,423){\vector(-1,0){80}}
              \put(398,503){\vector(0,-1){80}}
              \put(318,503){\vector(1,0){80}}
              \put(359,464){\vector(0,1){79}}
              \put(278,543){\vector(0,-1){80}}
              \put(359,544){\vector(-1,0){80}}
              \put(279,462){\vector(1,1){39}}
              \put(318,422){\vector(1,1){40}}
              \put(319,282){\vector(1,1){40}}
              \put(280,322){\vector(1,1){39}}
              \put(359,403){\vector(-1,0){80}}
              \put(279,403){\vector(0,-1){80}}
              \put(360,324){\vector(0,1){79}}
              \put(319,363){\vector(1,0){80}}
              \put(399,363){\vector(0,-1){80}}
              \put(399,283){\vector(-1,0){80}}
              \put(361,362){\circle{22}}
              \put(316,146){\vector(0,1){80}}
              \put(358,266){\vector(-1,0){80}}
              \put(278,266){\vector(0,-1){80}}
              \put(277,187){\vector(1,0){82}}
              \put(359,187){\vector(0,1){79}}
              \put(318,226){\vector(1,0){80}}
              \put(398,226){\vector(0,-1){80}}
              \put(398,146){\vector(-1,0){80}}
              \put(360,225){\circle{22}}
              \put(85,505){\circle{22}}
              \put(123,426){\vector(-1,0){80}}
              \put(123,506){\vector(0,-1){80}}
              \put(43,506){\vector(1,0){80}}
              \put(43,474){\vector(0,1){32}}
              \put(43,426){\vector(0,1){34}}
              \put(84,467){\vector(0,1){79}}
              \put(2,467){\vector(1,0){82}}
              \put(3,546){\vector(0,-1){80}}
              \put(83,546){\vector(-1,0){80}}
\end{picture}}

\begin{center}
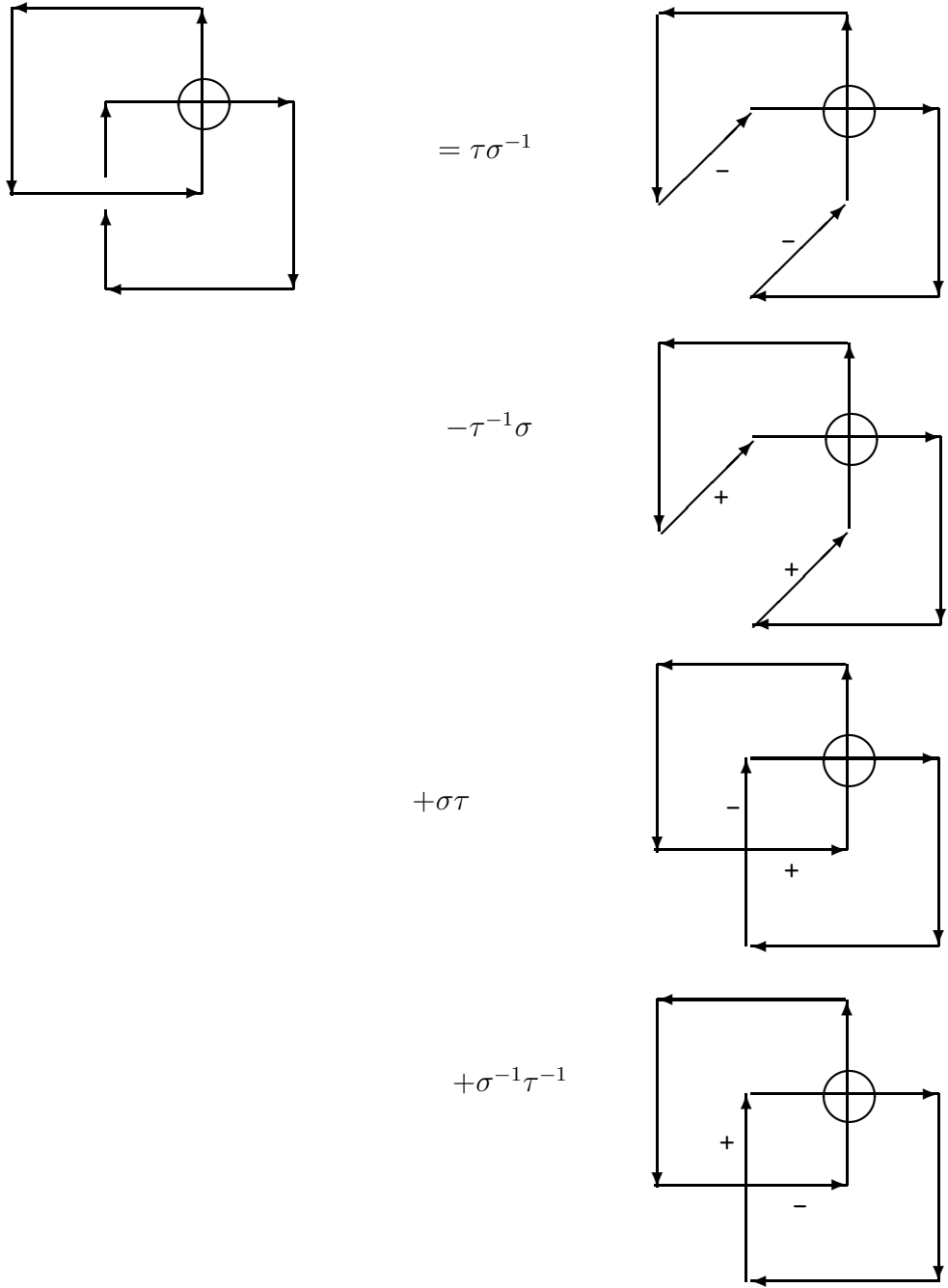
{\bf Figure 13  - Expansion Formulas for Virtual Hopf Link} \end{center}
\vspace{3mm}

\section{Bi-oriented Quantum Algebras} 
It is the purpose of this section to place our work with the generalized Alexander polynonmial in a context of bi-oriented quantum algebras.
To do this (and to define the concept of a bi-oriented quantum algebra) we need to first recall the notion of an oriented quantum algebra.
An {\em oriented quantum algebra} $(A, \rho, D, U)$ is an abstract
model for an oriented quantum invariant of classical links \cite{KRO}, \cite{KRCAT}. This model is based
on a solution to the Yang-Baxter equation and some extra structure
that serves to make an invariant possible to construct. The
definition of an oriented quantum algebra is as follows:   We are
given an algebra $A$ over a base ring $k$, an invertible solution
$\rho$ in $A \otimes A$ of the Yang-Baxter equation (in the
algebra formulation of this equation -- see the Remark below), and commuting automorphisms 
$U,D:A \longrightarrow A$  of the algebra, such that 

$$(U \otimes U)\rho = \rho,$$

$$(D \otimes D)\rho = \rho,$$
 
$$[(1_{A} \otimes U)\rho)][(D \otimes 1_{A^{op}})\rho^{-1}] 
= 1_{A \otimes A^{op}},$$ 

\noindent and 

$$[(D \otimes 1_{A^{op}})\rho^{-1}][(1_{A} \otimes U)\rho)]
= 1_{A \otimes A^{op}}.$$ 

\noindent The last two equations say that $[(1_{A} \otimes U)\rho)]$ and $[(D \otimes 1_{A^{op}})\rho^{-1}]$
are inverses in the algebra $A \otimes A^{op}$ where $A^{op}$ denotes the opposite algebra.  
\vspace{3mm}

When $U=D=T$, then $A$ is said to be {\em balanced}.  In this case
$$(T \otimes T)\rho = \rho,$$

$$[(1_{A} \otimes T)\rho)][(T \otimes 1_{A^{op}})\rho^{-1}] = 1_{A
\otimes A^{op}}$$

\noindent and 

$$[(T \otimes 1_{A^{op}})\rho^{-1}][(1_{A} \otimes T)\rho)] = 1_{A
\otimes A^{op}}.$$

In the case where $D$ is the identity mapping, we call the 
oriented quantum algebra {\em standard}.  As we saw in  \cite{KRCAT}, 
the invariants defined by Reshetikhin and Turaev (associated with a 
quasi-triangular Hopf algebra) arise from standard oriented quantum algebras.
It is an interesting structural feature of algebras that we have elsewhere 
\cite{GAUSS} called {\em quantum algebras} (generalizations of quasi-triangular Hopf algebras)
that they give rise to standard oriented quantum algebras. Note that the term quantum algebra
as used here is more specific than the $QA$ quantum algebra designation that is used in archived 
papers related to Hopf algebras and quantum groups. 
\vspace{3mm}

Appropriate matrix representations of oriented quantum 
algebras or the existence of certain traces on these algebras
allow the construction of oriented invariants of knots and links.
These invariants include all the known quantum link
invariants at the time of this writing. \vspace{3mm}

\noindent {\bf Remark.} Note that we have the Yang-Baxter elements $\rho$ and $\rho^{-1}$
in $A \otimes A.$ We assume that $\rho$ and $\rho^{-1}$ satisfy the algebraic Yang-Baxter
equation. This equation (for $\rho$) states

$$\rho_{12}\rho_{13}\rho_{23} = \rho_{23}\rho_{13}\rho_{12}$$

\noindent where $\rho_{ij}$ denotes the placement of the tensor factors of $\rho$ in the 
$i$-th and $j$-th tensor factors of the triple tensor product $A \otimes A \otimes A.$
\vspace{3mm}

\noindent We write $\rho = \Sigma e \otimes e'$ and
$\rho^{-1} = \Sigma E \otimes E'$ to indicate that these elements are sums of
tensor products of elements of $A$. The expression $e \otimes e'$ is thus a generic 
element of the tensor product. However, we often abbreviate and write
$\rho = e \otimes e'$ and $\rho^{-1} = E \otimes E'$ where the summation is implicit. 
We refer to $e$ and $e'$ as the {\em signifiers} of $\rho$,
and $E$ and $E'$ as the {\em signifiers} of $\rho^{-1}$. For example, 
$\rho_{13} = e \otimes 1 \otimes e'$ in $A \otimes A \otimes A.$ 
\vspace{3mm}

\noindent Braiding operators, as they appear in knot theory, differ from the algebraic 
Yang-Baxter elements by a permutation of tensor factors. This point is crucial to the 
relationship of oriented quantum algebras and invariants of knots and links.
\vspace{3mm} 

\noindent We extend the concept of oriented quantum algebra  by adding a second solution to 
the Yang-Baxter equation $\gamma$ that will take the role of the virtual crossing.
\smallbreak

\noindent {\bf Definition.} A {\em bi-oriented quantum algebra} is a quintuple  $(A, \rho, \gamma, D , U)$ such that
$(A, \rho, D, U)$ and $(A,\gamma, D, U)$ are oriented quantum algebras and $\gamma$ satisfies the following properties:

\begin{enumerate}
\item $\gamma_{12}\gamma_{21} = 1_{A \otimes A}.$ (This is the equivalent to the statement that the 
braiding operator corresponding to $\gamma$ is its own inverse.) 

\item The following mixed identities involving $\rho$ and $\gamma$ are satisfied. These correspond to the 
braiding versions of the virtual detour move of type three that involves two virtual crossings and one
real crossing.

$$\gamma_{12}\rho_{13}\gamma_{23} = \gamma_{23}\rho_{13}\gamma_{12}$$

$$\gamma_{12}\gamma_{13}\rho_{23} = \rho_{23}\gamma_{13}\gamma_{12}$$

$$\rho_{12}\gamma_{13}\gamma_{23} = \gamma_{23}\gamma_{13}\rho_{12}.$$
 
\end{enumerate}
\bigbreak

By extending the methods of \cite{KRCAT}, it is not hard to see that {\em a bi-oriented quantum algebra will always give rise to invariants
of virtual links up to the type one moves (framing and virtual framing).} 
\bigbreak 

\noindent Any oriented quantum algebra $(A, \rho, D, U)$ gives rise to a 
bi-oriented quantum algebra by taking $\gamma$ to be the identity element in $A \otimes A.$ This corresponds 
to associating a simple permutation (transposition) to the virtual crossing. The resulting invariants of virtuals are worth investigating. In 
\cite{VKT} the corresponding generalization of the Jones polynomial is studied. The bi-oriented quantum algebras form a context for these invariants. 
\bigbreak

In the case of the generalized Alexander polynomial studied
in this paper, the matrix model and state model for $Z(K)$ translate directly into a specific example of a bi-oriented balanced quantum
algebra $(A, \rho, \gamma, T)$ (It is a balanced bi-oriented quantum algebra.) with the underlying algebra $A$ 
the algebra of elementary matrices as in \cite{KRCAT}. In making this translation to the 
algebra, one must take the matrices of the matrix model and compose with a permutation. In our case the matrix $\gamma$ is a 
diagonal $4 \times 4$ matrix with three ones and one minus one in that order on the diagonal. This is the matrix obtained from the matrix
$\eta$ that we used for the virtual crossing.
\bigbreak

$$\gamma =  \left[
\begin{array}{cccc}
     1 & 0 & 0 & 0 \\
     0 & 1 & 0 & 0\\
     0 & 0 & 1 & 0 \\
     0 & 0 & 0 & -1
\end{array}
\right]. $$ 

\noindent The matrix $\rho$ is obtained from the braiding matrix $R$ by permuting the two middle columns (because $\rho^{ab}_{cd} = R^{ba}_{cd}$).

$$\rho =  \left[
\begin{array}{cccc}
     \sigma^{-1} \tau & 0 & 0 & 0 \\
     0 & \sigma \tau & z & 0\\
     0 & 0 & \sigma^{-1} \tau^{-1} & 0 \\
     0 & 0 & 0 & -\sigma \tau^{-1}
\end{array}
\right]. $$

\noindent On elementary matrices $E^{a}_{b}$ the transformation $T$ is given by the formula $T(E^{a}_{b}) = i^{b-a}$ where $i^{2} = -1.$
(See Section $4.2$ of \cite{KRCAT}.)  This completes the description of the bi-oriented quantum algebra that corresponds to the generalized Alexander
polynomial of this paper.   
\bigbreak

The main algebraic point about the bi-oriented quantum algebra for the generalized Alexander polynomial is that the operator $\gamma$ for the virtual
crossing is not the identity operator, and that this non-triviality is crucial to the structure of the invariant.
We will investigate bi-oriented quantum  algebras and other examples of virtual invariants derived from them in a subsequent paper.
\bigbreak

 \end{document}